\documentclass{article}
\usepackage{fullpage}
\usepackage{amsmath}
\usepackage{amsfonts}
\usepackage{epsfig}
\usepackage{graphics} 

\usepackage{psfrag}
\usepackage{color}
\usepackage{amssymb}

\usepackage{cite}

\newcommand{\R}{\mathbb R}

\newcommand{\N}{\mathbb N}

\renewcommand{\P}{\mathbb P}

\newcommand{\dpar}[2]{\dfrac{\partial #1}{\partial #2}}

\newcommand{\bu}{\mathbf{u}}

\newcommand{\bF}{\mathbf{f}}

\newcommand{\bn}{\mathbf{n}}
\newcommand{\bx}{\mathbf{x}}

\newcommand{\bug}[1]{{#1}}

\begin{document}

\title{How to avoid mass matrix for linear hyperbolic problems}


\author{ R\'emi Abgrall\footnote{remi.abgrall@math.uzh.ch}, Paola Bacigaluppi and  Svetlana Tokareva \\
  Institute of Mathematics, University of Z\"urich, \\Winterthurerstrasse 190, CH 8057 Z\"urich, Switzerland,  
  }

\maketitle

\begin{abstract}
  We are interested in the numerical solution of linear hyperbolic problems using continuous finite elements of arbitrary
  order. It is well known that this kind of methods, once the weak formulation has been written, leads to a system of   ordinary differential equations in $\R^N$, where $N$ is the number of degrees of freedom. The solution of the resulting ODE system involves the inversion of a
  sparse mass matrix that is not block diagonal. Here we show how to avoid this step, and what are the
  consequences of the choice of the finite element space. Numerical examples show the correctness of our approach.
\end{abstract}

\section{Introduction}
We are interested in the numerical approximation of the hyperbolic problem
\begin{subequations}
\label{abgrall_pleinary_convection}
\begin{equation}\label{abgrall_pleinary_eq:1}
\dpar{u}{t}+\text{ div }\bF(\bx, u)=0 \qquad \bx\in\Omega\subset \R^d
\end{equation}
by means of a finite element like technique. In this paper,
we focus on the linear case where $\bF(\bx,u)=\mathbf{a}(\bx)u$. The vector field $\mathbf{a}$ may depend on the
spatial location $\bx$. The problem \eqref{abgrall_pleinary_eq:1} is also supplemented with initial and boundary conditions:
\begin{equation}
\label{abgrall_pleinary_eq:1:2}
u(\bx, 0)=u_0(\bx)
\end{equation}
and 
\begin{equation}
\label{abgrall_pleinary_eq:1:3}
u(\bx,0)=g(\bx) \text{ if } \bx \in \partial \Omega.
\end{equation}
Obviously, \eqref{abgrall_pleinary_eq:1:3} has to be understood in the weak sense, i.e. that $u=g$ on the inflow characteristics.
\end{subequations} 

The physical space is covered by a conformal tessellation $\mathcal{T}$. For ease of exposition, we assume that
$$\Omega=\cup_{K\in \mathcal{T}} K.$$
 The solution of the problem is approximated by an element of the space $V^h$ defined by:
$$V^h=\{ \bu^h\in C^0(\Omega) \text{ such that for any }K, \bu^h_{|K}\text{ is a polynomial of degree }r\}.$$
We denote by $\mathbb{P}^r$ the set of polynomials of degree $r$. In this paper, we consider $r=1,2$ only.

It is well known that any finite element technique applied to \eqref{abgrall_pleinary_eq:1} will lead to a formulation of the type
$$M \dfrac{dU}{dt}+F=0$$
where $U$ denotes the vector of degrees of freedom, $F$ is an approximation of the term $\mathrm{div}\,\bF$ and $M$ is a mass
matrix. In the case of continuous elements, this matrix is sparse but not block diagonal, contrarily to what happens
for the Discontinuous Galerkin methods where the global continuity requirement is not made. Hence, in order to use any
standard ODE solver, we need to invert $M$. This is considered cumbersome by many practitioners and this has been,
in our opinion, one of the factors that has led to supremacy of DG methods in the current development of high order
schemes.

Several researchers have proposed methods that avoid this step. More precisely, their methods are designed in such a way that the actual mass matrix is diagonal, so that the problem amounts to finding a ``good'' lumping integration formula. 
The first work we are aware of in that direction is \cite{abgrall_pleinary_ho_lumping1}, where the wave equation is considered, and the finite
element space is made of functions belonging to a \bug{subspace of $\P^{k+1}$ that contains $\P^k$}. This amounts to adding one degree of freedom to the ``natural'' quadratic elements. This work has been followed,
in the same spirit, by \cite{abgrall_pleinary_ho_lumping2} where higher accuracy could be  obtained. However, the elements become more and more complex and, what is even more important, the stability condition on the time step becomes dramatically restrictive.

In these notes we describe some preliminary results about a new method for which no inversion of the mass
matrix is needed, while  a typical finite element  approximation can be kept for the description of the divergence term. In this approach, there is no need to change the degrees of freedom. The method presented here can be seen as an extension of \cite{abgrall_pleinary_mario} where only $\P^1$ elements and second order approximation in time have been considered.

The rest of the paper is organized as follows. In the first section, we describe the approximation of the divergence term of \eqref{abgrall_pleinary_eq:1}. These are classical 
stabilized finite element methods. In the second section, we describe and somewhat justify our approach. The last section provides numerical examples that justify the correctness of our approach. A more involved analysis and and description will be made elsewhere. We conclude by giving some perspectives.

 \section{Description of the scheme}
 We start by describing the two spatial approximations  we     consider, then explain how to avoid the mass matrix inversion.
We are given a triangulation of $\R^d$. Here we assume $d=2$, but the discussion is general. The elements are denoted
by $K$ and assumed to be simplices. In each element, we assume that the solution is approximated by a
polynomial of degree $r$ and that the approximation is globally continuous. Let us denote the approximate solution by $u^h$.
The function $u^h$ is fully defined by its control parameter $u_\sigma$ at all the degrees of freedom $\sigma$.
We define by $\mathcal{S}$ the set of degrees of freedom, so that
    $$u^h=\sum_{\sigma\in \mathcal{S}} u_\sigma \varphi_\sigma.$$
    
We denote by $V_h=\text{span }\{\varphi_\sigma, \sigma\in \mathcal{S}\}$
For now, we can think of $u_\sigma$ as the value of $u^h$ at $\sigma$ and thus $\varphi_\sigma$ is the Lagrange basis, but we will need slightly less conventional approximation later.

    We assume that we have a good integrator of the steady version of \eqref{abgrall_pleinary_convection}, and that this scheme writes: for any degree of freedom $\sigma$, $u^h$ satisfies:
    $$\sum_{K\ni \sigma} \Phi_\sigma^{K,\bx}(u^h)=0.$$

    Examples are given by:
    \begin{enumerate}
    \item The SUPG residual, \cite{abgrall_pleinary_JohnsonNavertPitaka,abgrall_pleinary_Hughes1987}:
      \begin{equation}\label{abgrall_pleinary_SUPG}
        \begin{split}
          \Phi_\sigma^\bx(u^h)&=\int_{\partial K}\varphi_\sigma \bF(u^h)\cdot \bn \;d\mathbf{\ell} -\int_K \nabla \varphi_\sigma\cdot \bF(u^h)\; d\bx \\
          &\qquad +h_K
          \int_K \bigg (\nabla_u\bF(u^h)\cdot \nabla \varphi_\sigma \bigg )\tau \bigg (\nabla_u\bF(u^h)\cdot \nabla u^h \bigg )\; d\bx
        \end{split}
        \end{equation}
    with $\tau>0$.  We take:
$$(h_K \tau)^{-1}=\sum\limits_{\sigma\in K} |\overline{\mathbf{a}}_K\cdot \nabla\varphi_\sigma|$$ where
$\overline{\mathbf{a}}_K$ is the value of $\mathbf{a}$ at the centroid of $K$.
    \item The Galerkin scheme with jump stabilization \cite{abgrall_pleinary_burman}:
      \begin{equation}\label{abgrall_pleinary_Jump}
        \begin{split}
          \Phi_\sigma^\bx(u^h)&=\int_{\partial K}\varphi_\sigma \bF(u^h)\cdot \bn \;d\mathbf{\ell} -\int_K \nabla \varphi_\sigma\cdot \bF(u^h) \; d\bx \\
          &\qquad +\sum_{edges}
          \Gamma h_e^2 \int_e [\nabla u]\cdot [\nabla \varphi_\sigma]^+\;d\mathbf{\ell}
        \end{split}
        \end{equation}
    with $\Gamma>0$. 
    Here, since the mesh is conformal, any edge (or face in 3D) is the intersection of the element $K$ and an other element denoted by $K^+$.
    We define $[\nabla u]=\nabla u_{|K}-\nabla u_{| K^+}$ and $[\nabla \varphi_\sigma]^+=(\varphi_\sigma)_{|K}$.
Here, we have taken $\Gamma=\max(\overline{\mathbf{a}}_K, \overline{\mathbf{a}}_{K^+} )$.
 See \cite{abgrall_pleinary_burman} for more details.
    \end{enumerate}
    This streamline formulation implies formally that the exact solution cancels the  residuals.
    In the case of the stabilisation by jumps, we can only write that
    $$\Phi_\sigma^K=\int_K\varphi_\sigma \text{div }\bF(u) dx+R_\sigma(u^h)$$
    where
    $\sum_{\sigma \in K} R_\sigma(u^h)=0$.
    The additional term $R_\sigma$ is non-zero, except for the exact solution unless this solution has continuous normal gradients. 
     For steady solutions, both methods can be shown to converge as $h^{k+1/2}$, see \cite{abgrall_pleinary_JohnsonNavertPitaka,abgrall_pleinary_burman} for more details.

    \subsection{Formulation for unsteady problems}
     We use a deferred correction (DeC) approach. We start from the ODE:
    \begin{equation}
    \label{DeC:ODE}    
    \dfrac{dy}{dt}=f(y,t), \qquad y(0)=y_0.
    \end{equation}
    We follow the main ideas of \cite{abgrall_pleinary_greengard}. Between $t_n$ and $t_{n+1}$, the solution of \eqref{DeC:ODE} satisfies
    $$y(t)=y(t_n)+\int_{t_n}^t f(y(s),s) ds.$$
    Given $0=\xi_0<\xi_1<\ldots <\xi_l<\ldots <\xi_{M+1}=1$, and we consider the times $t_{n,l}=t_n+\xi_l \Delta t$ with $\Delta t=t_{n+1}-t_n$.
    If we know $f_{n,l}\approx f(y(t_{n,l}),t_{n,l})$, we can consider the Lagrange interpolant $\mathcal{I}_{M+1}$ of $f$ with data given by 
$(t_{n,l}, f_{n,l})$, therefore we get the approximation:
    $$y_{n,l}=y_{n,0}+\int_{t_n}^{t_{n,l}} \mathcal{I}_{M+1}[f(y(~.~),~.~)](t_n+\xi\Delta t) d\xi.$$
    This is in general a non-linear implicit equation. 
    
    The idea of the DeC method is to consider the first order scheme, for
    $M\geq l\geq 1$:
    $$y_{n,l}=y_{n,l-1}+ \alpha_l\Delta t f(y(t_{n,l-1}),t_{n,l-1}), \qquad y_{n,0}\approx y(t_n)$$
    where $\alpha_l=\xi_{l}-\xi_{l-1}$.
    Then, we introduce the vector $v=(y_{n,1}, \ldots , y_{n,M+1})^T$. The first order scheme can be rewritten as
    $L^1(v)=0$ where
    $$L^1(v)=\begin{pmatrix}
    y_{n,1}-y_{n,0}-\Delta t \int_0^{\xi_1}\mathcal{I}_0[f(y(~.~),~.~)](t_n+\xi \Delta t) d\xi\\
    \vdots \\
    y_{n,l}-y_{n,0}- \Delta t \int_0^{\xi_{l}}\mathcal{I}_0[f(y(~.~),~.~)](t_n+\xi \Delta t) d\xi\\
    \vdots \\
    y_{n,M+1}-y_{n,0}- \Delta t \int_0^{\xi_M}\mathcal{I}_0[f(y(~.~),~.~)](t_n+\xi \Delta t) d\xi
    \end{pmatrix}
    $$
    where $\mathcal{I}_0$ is the first order interpolant of $f$: for $1 \leq l \leq M+1$, 
    {$$\mathcal{I}_0[f(y(~.~),~.~)](s)=f(y_{n,l-1},t_{n,l-1}) \qquad \text{ for } s\in [t_{n,l-1},t_{n,l}[.$$}
Note that $L^1(v)=0$ can be solved \emph{explicitely}.

Similarly, we define $L^2$ by:
    $$L^2(v)=\begin{pmatrix}
    y_{n,1}-y_{n,0}-\Delta t \int_0^{\xi_1}\mathcal{I}_{M+1}[f(y(~.~),~.~)](t_n+\xi \Delta t) d\xi\\
    \vdots \\
    y_{n,l}-y_{n,0}- \Delta t \int_0^{\xi_{l}}\mathcal{I}_{M+1}[f(y(~.~),~.~)](t_n+\xi \Delta t) d\xi\\
    \vdots \\
    y_{n,M+1}-y_{n,0}- \Delta t \int_0^{\xi_M}\mathcal{I}_{M+1}[f(y(~.~),~.~)](t_n+\xi \Delta t) d\xi
    \end{pmatrix}.
    $$
 
          The DeC formulation is defined as follows:
          \begin{enumerate}
          \item $v^0=(y_{n}, \ldots y_n)^T$ and $y_{n,0}=y_n$,
          \item For $ k=1, \ldots M+1$, $v^{k}$ is defined as
            $$L^1(v^{k})=L^1(v^{k-1})-{L^2(v^{k-1})}$$
          \end{enumerate}
          Since $L^1$ is explicit, the method is completely explicit. One can show that $L^2-L^1=O(\Delta t)$ so that the scheme is $(M+1)$-th order accurate.

          \bigskip 
          
    Similar to what is done for  ODEs, we could integrate \eqref{abgrall_pleinary_convection} in time and get:
    $$u(\bx,t_{n+1})=u(\bx,t_n)+\int_{t_n}^{t_{n+1}} \text{div }\bF(u(x,s)) ds,$$
    This can be approximated by 
    \begin{equation}
    \label{abgrall_pleinary_ODE}
\begin{split}
    u(\bx,t_n+\xi_i\Delta t)&\approx u(\bx,t_n)+\int_0^{\xi_i} \text{div }\mathcal{I}_{r+1}[\bF(u(\bx, ~.~) ](t_n+\xi \Delta t) ds \\&\qquad = \Delta t\;\sum_{l=0}^{r} \omega_l^i \text{div }\bF(u(\bx, \xi_j) ds
\end{split}
    \end{equation}
    $\mathcal{I}_{r+1}[\bF(u(\bx, ~.~)) ]$ is the Lagrange interpolant of $\bF(u(\bx, ~.~))$  at the points $\{t_n, , \ldots,  \xi_i \Delta t, \ldots,  t_{n+1}\}$ and 
     $\omega_l^i$ are the weights.
     
      This suggests the algorithm we describe now. For any $V\in V_h^M$, $V^{\sigma}=(V_1^\sigma, \ldots , V_{M+1}^\sigma)^T$ is a
 vector of  control
    parameters at the degree of freedom $\sigma\in \mathcal{S}$: $V=\sum_{\sigma\in \mathcal{S}} V^\sigma\varphi_\sigma$.
        Then,  we can consider the following deferred correction approximation: we introduce $t_{n,i}=t_n+\xi_i (t_{n+1}-t_n)$ so that $t_{n,0}=t_n$ and $t_{n,r+1}=t_{n+1}$, and define 
 \begin{enumerate}
 \item for any $\sigma\in \mathcal{S}$, the operator $L^1_\sigma$ as
 \begin{subequations}\label{abgrall_pleinary_2}
 \begin{equation}
 \label{abgrall_pleinary_2.1}
 L^1_\sigma(V_1, \ldots, V_{r+1})=\begin{pmatrix}
 |C_\sigma|(V_{r+1}^\sigma-V_0^\sigma)+\sum\limits_{K\ni\sigma}\displaystyle\int_{t_{n,0}}^{t_{n,r+1}}  \mathcal{I}_0[\Phi_\sigma^\bx ] \big(t_n+s\Delta t \big) \;ds\\
 |C_\sigma|(V_{r}^\sigma-V_{0}^\sigma)+\sum\limits_{K\ni\sigma}\displaystyle\int_{t_{n,0}}^{t_{n,r}} \mathcal{I}_0[\Phi_\sigma^\bx ] \big(t_n+s\Delta t)\big) \;ds\\
\vdots\\
 |C_\sigma|(V_{1}^\sigma-V_0^\sigma)+\sum\limits_{K\ni\sigma}\displaystyle\int_{t_{n,0}}^{t_{n,1}} \mathcal{I}_0[\Phi_\sigma^\bx] \big(t_n+s\Delta t\big) \;ds\\
 \end{pmatrix}
\end{equation}
Here, $V_0^\sigma=(u^n_\sigma, \ldots , u^n|\sigma)^T\in \R^M$.
\item and the operator $L^2_\sigma$ as
\begin{equation}
\label{abgrall_pleinary_2.2}L^2_\sigma(V_1, \ldots, V_{r+1})=\begin{pmatrix}
\sum\limits_{K\ni\sigma}\Bigg (\displaystyle\int_K \Psi_\sigma \big ( V_{r+1}-V_0\big ) \; dx+\displaystyle\int_{t_{n,0}}^{t_{n,r+1}} \mathcal{I}_{r+1}[\Phi_\sigma^\bx ] \big(t_n+s\Delta t \big) \; ds\\

\sum\limits_{K\ni\sigma}\Bigg (\displaystyle\int_K \Psi_\sigma \big ( V_{r}-V_0\big )  \;dx+\displaystyle\int_{t_{n,0}}^{t_{n,r}}  \mathcal{I}_{r+1}[\Phi_\sigma^\bx ] \big(t_n+s\Delta t \big)\; ds\\

\vdots\\
\sum\limits_{K\ni\sigma}\Bigg (\displaystyle\int_K \Psi_\sigma \big ( V_{1}-V_0\big ) \; dx+\displaystyle\int_{t_{n,0}}^{t_{n,1}}  \mathcal{I}_{r+1}[\Phi_\sigma^\bx ] \big(t_n+s\Delta t \big)\; ds
\end{pmatrix}
\end{equation}
\end{subequations}
 \end{enumerate}

Last, we define the operators $L^1$ and $L^2$ on the finite element set $V_h$ as
$$L^1=(L^1_\sigma)_{\sigma\in \mathcal{S}}, \qquad L^2=(L^1_\sigma)_{\sigma\in \mathcal{S}}.$$

The step of the method between $t_n$ and $t_{n+1}$ is defined as follows.
\begin{enumerate}
\item Knowing $u_\sigma^n$, we set $V^0_\sigma=(u^n_\sigma, \ldots , u^n_\sigma)$.
\item For $k=1,\dots,M$, we construct $V^k$ as the solution of
$$L^1(V^{k+1})=L^1(V^k)-L^2(V^k).$$
\item Then we define $u_\sigma^{n+1}$ as
$$u_\sigma^{n+1}=(V_\sigma^{r+1})^M.$$
\end{enumerate}

This methods provides a decent approximation of the solution because one can show \cite{abgrall_pleinary_abgrall:high} that, for the $L^2$ norm,
\begin{equation}
\label{abgrall_pleinary_condition}
||L^1-L^2||\leq C \Delta t,
\end{equation}
where the constant $C$ depends only on the mesh.
Then, using standard results for deferred correction methods, one can show that we have an $(r+1)$-th order accurate scheme if $M=r+1$, provided $L^1$ is invertible. The overall cost is not larger than a standard Runge-Kutta method.

Let us now have a look at the invertibility of $L^1$.
Not every finite element approximation can work. The reason is that we have not yet specified what should be the parameters
 $C_\sigma$ in relation \eqref{abgrall_pleinary_2.1}. It is easy to see that we must have
$$C_\sigma=\int_\Omega \varphi_\sigma d\bx,$$
and in order that $L^1$ be invertible, we need $C_\sigma>0$. For $\P^1$ elements, there is no problem because  the basis functions are positive, but it is well known that this condition is not met for higher order finite elements. For example, in the case of two-dimensional quadratic Lagrange interpolation, we have six basis functions.
 Three of them are associated to the vertices, and it is well known that their integral vanishes, so that in the end $C_\sigma=0$ for the vertices. For other finite elements, we can have
$C_\sigma<0$.

In order to circumvent this restriction, and since we are interested in the approximation order and \emph{not on the practical representation, i.e. the physical meaning
 of the degrees of freedom}, a simple  way is to replace classical Lagrange elements of degree $r$ by their Bezier counterparts.
If 
$$\bigg (\sum\limits_{j=1}^{d+1}x_j\bigg )^{r}=\sum\limits_{\sum_{k=1}^{d+1} j_k=r}\theta_{j_1\ldots j_{d+1} }^r x_1^{j_1}\ldots x_{j_{d+1}}$$
is the binomial expansion, then the Bezier polynomials are simply
$$B_{j_1\ldots j_{d+1}}^r=\theta_{j_1\ldots j_{d+1}}^r \Lambda_1^{j_1}\ldots \Lambda_{j_{d+1}}$$
where the $\Lambda_j$ are the standard barycentric coordinates. Since 
$$\int_KB_{j_1\ldots j_{d+1}}^r(\bx) d\bx >0,$$
and since this family is a basis of $\P^r$, there are no more problems.  In the simulations done in this paper, we have chosen $\P^1$ elements (i.e. B\'ezier of degree 1), 
and quadratic B\'ezier elements. Note that this kind of approximation has already been used for steady problems \cite{abgrall_pleinary_jirka}, and has some links with 
isogeometrical analysis \cite{abgrall_pleinary_isoHughes}, but for completely different reasons.
\section{Numerical illustrations}
\subsection{Parameters}
In the numerical experiments we present, we have chosen a temporal scheme that is third order in time. It is based on the Lagrange interpolation in $[0,1]$, where the data are given at the points $t=0$, $\frac{1}{2}$ and $1$.
This results in the following formula that defines the operator $L^2$:
\begin{equation*}
\begin{split}
\int_0^{1/2} \mathcal{I}_2(f)ds&=\frac{5}{24}f(0)+\frac{1}{3}f(\frac{1}{2})-\frac{1}{24} f(1)\\
\int_0^1\mathcal{I}_2(f) ds&= \frac{1}{6}f(0)+\frac{4}{6}f(\frac{1}{2})+\frac{1}{6}f(1)
\end{split}
\end{equation*}
We have used the same temporal scheme for $\P^1$ and $\mathbb{B}^2$ elements.

\subsection{Simulations}
\begin{figure}[h]
\begin{center}
\begin{tabular}{cc}
\includegraphics[width=0.45\textwidth]{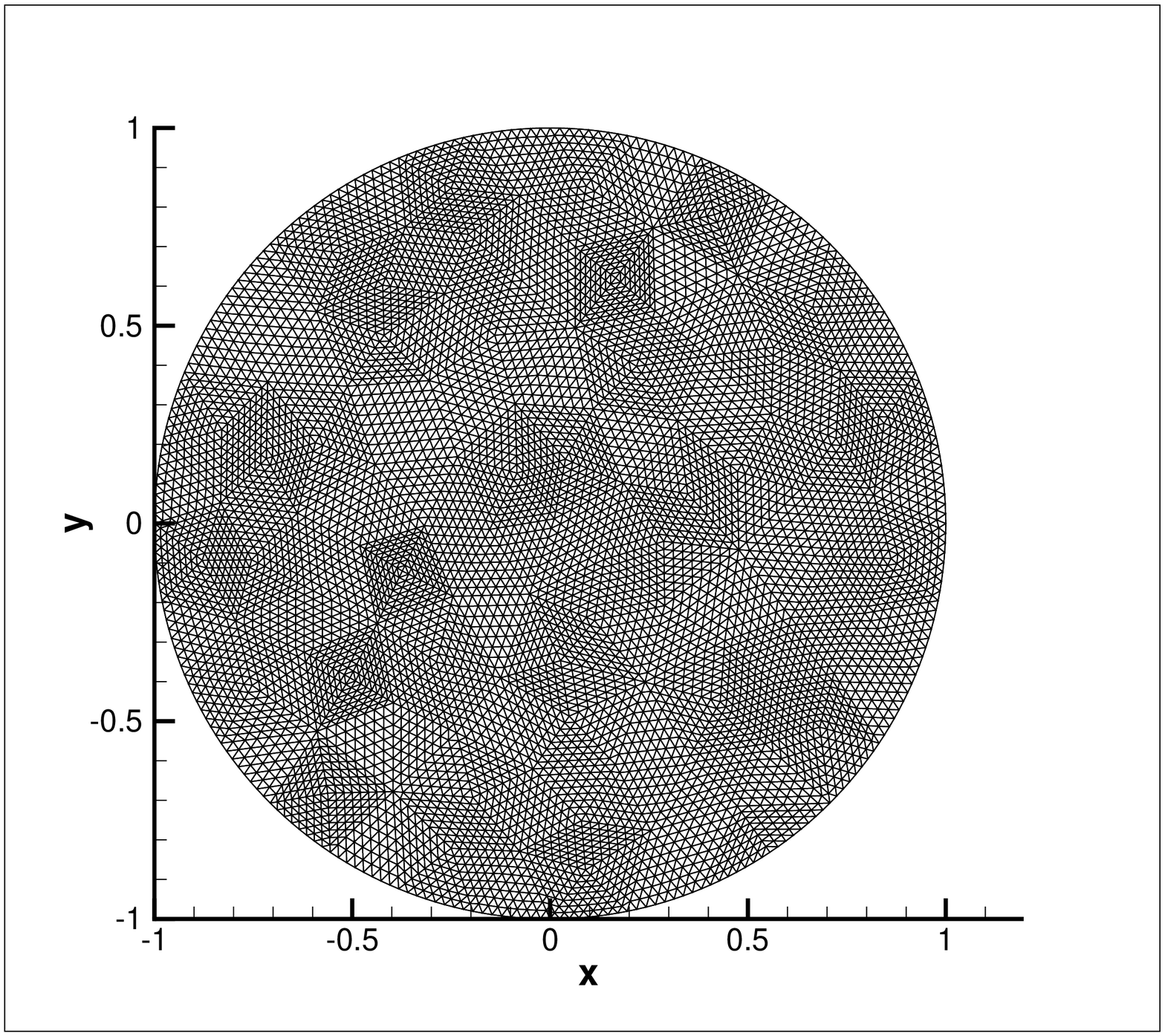}&
\includegraphics[width=0.45\textwidth]{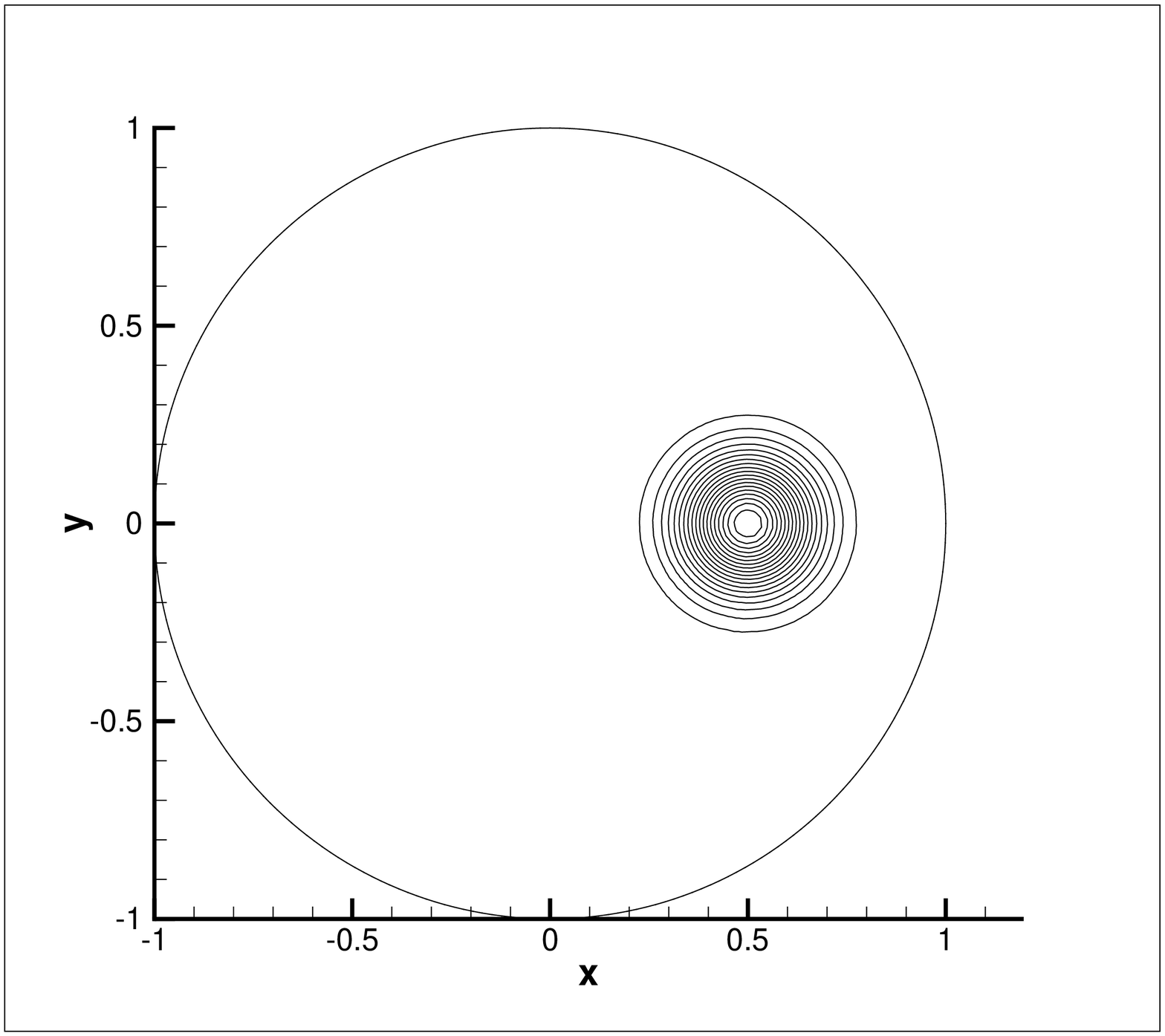}
\end{tabular}
\caption{\label{abgrall_pleinary_mesh}Exact solution after $n$ rotations ($n\in \N$) and plot of the degrees of freedoms.}
\end{center}
\end{figure}

The velocity field at $(x,y)$ is given by $\mathbf{a}=2\pi (-y,x)$. The initial condition is given by:
$$
u_0(x,y)=e^{-40(x^2+y^2)}
.
$$
The domain is a circle with center $(0,0)$ and radius $R=1$. The mesh representing all the degrees of freedom is displayed in Figure \ref{abgrall_pleinary_mesh}:
 The quadratic elements have $6$ degrees of freedom (the vertices and the mid-points of the edges).  \bug{These degrees of freedom are also used for
 the linear element just by mesh refinement}. There are $7047$ degrees of freedom here, so $h\approx\sqrt{\frac{\pi}{7047}}\approx 0.021$ which is
 relatively coarse. On the same figure, we represent the exact solution.
The time step is evaluated as the minimum of the $\Delta t_K$ defined by:
$$\Delta t_K=\text{CFL }\dfrac{h_K}{||\overline{\mathbf{a}}_K||}$$
where $h_K$ is the length of the smallest edge of $K$ and $\overline{\mathbf{a}}_K$ is the speed at the centroid. Since the elements for the $\P^1$ simulations
 are obtained from those of the $\mathbb{B}^2$ simulation by splitting, the parameter $h_K$, for the $\P^1$ simulations, is half of the one for the $\mathbf{B}^2$ 
simulations. \bug{For that reason, the CFL number for the quadratic approximation is half of the one chosen for the linear simulations, namely $0.6$ instead of $0.3$:
 we run with the same time step.} By the way, we have not yet conducted a rigorous study of the CFL condition, but all experiments
 indicate that the quadratic simulations can be safely run with $CFL=0.5$.

\begin{figure}[h]
\begin{center}
\begin{tabular}{cc}
\includegraphics[width=0.45\textwidth]{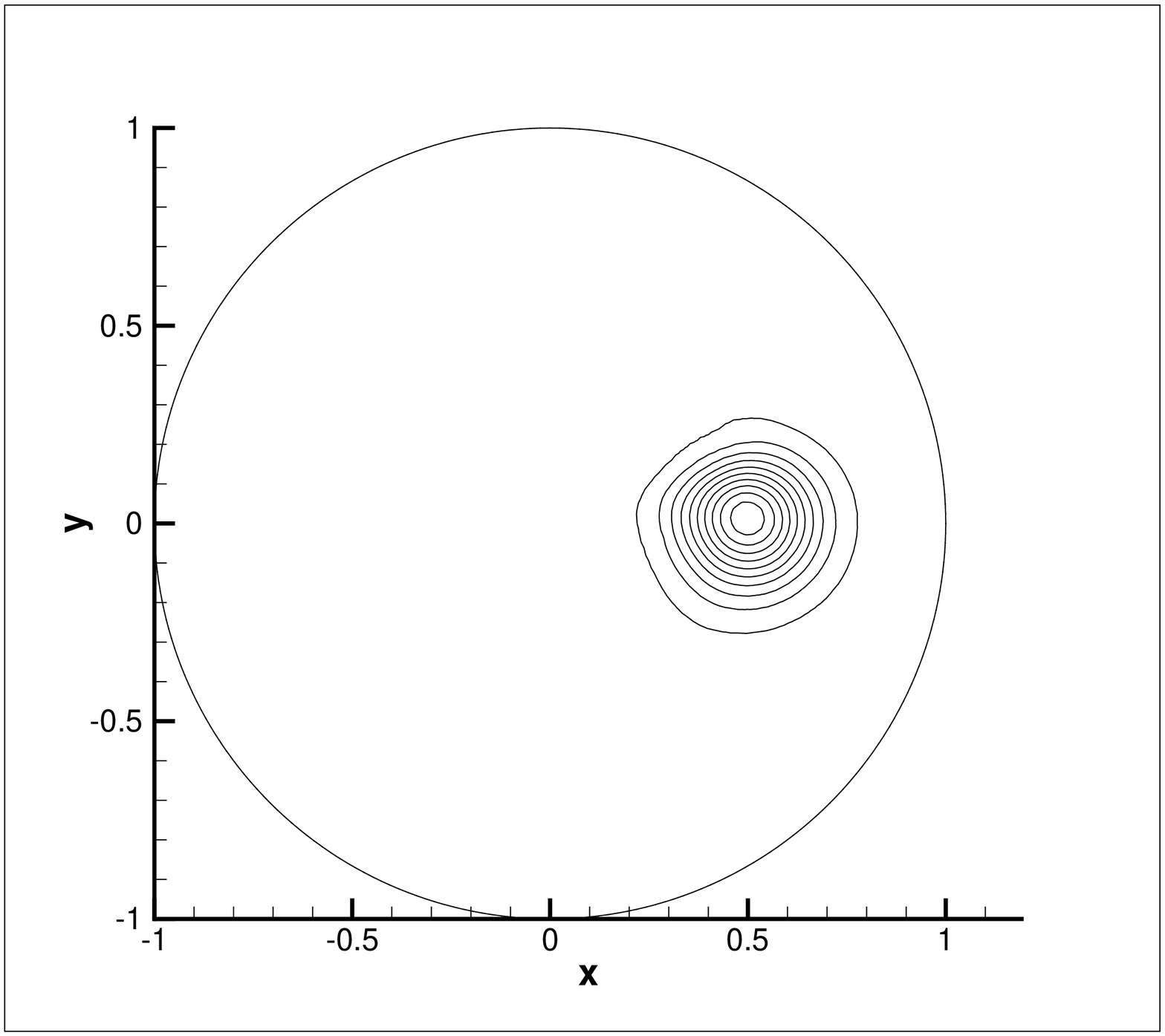}&
\includegraphics[width=0.45\textwidth]{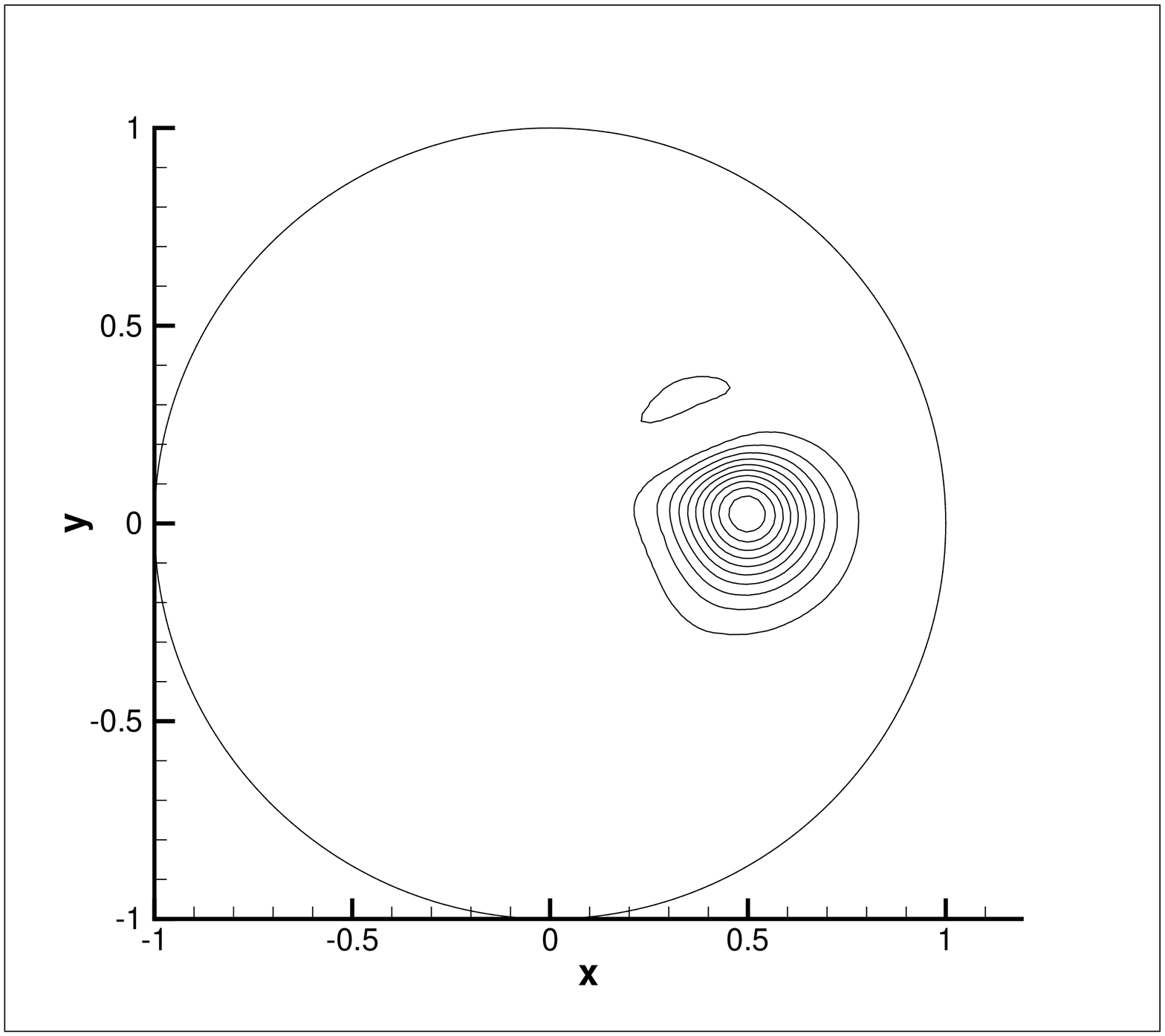}\\
(a)&(b)\\
\multicolumn{2}{c}{\includegraphics[width=0.45\textwidth]{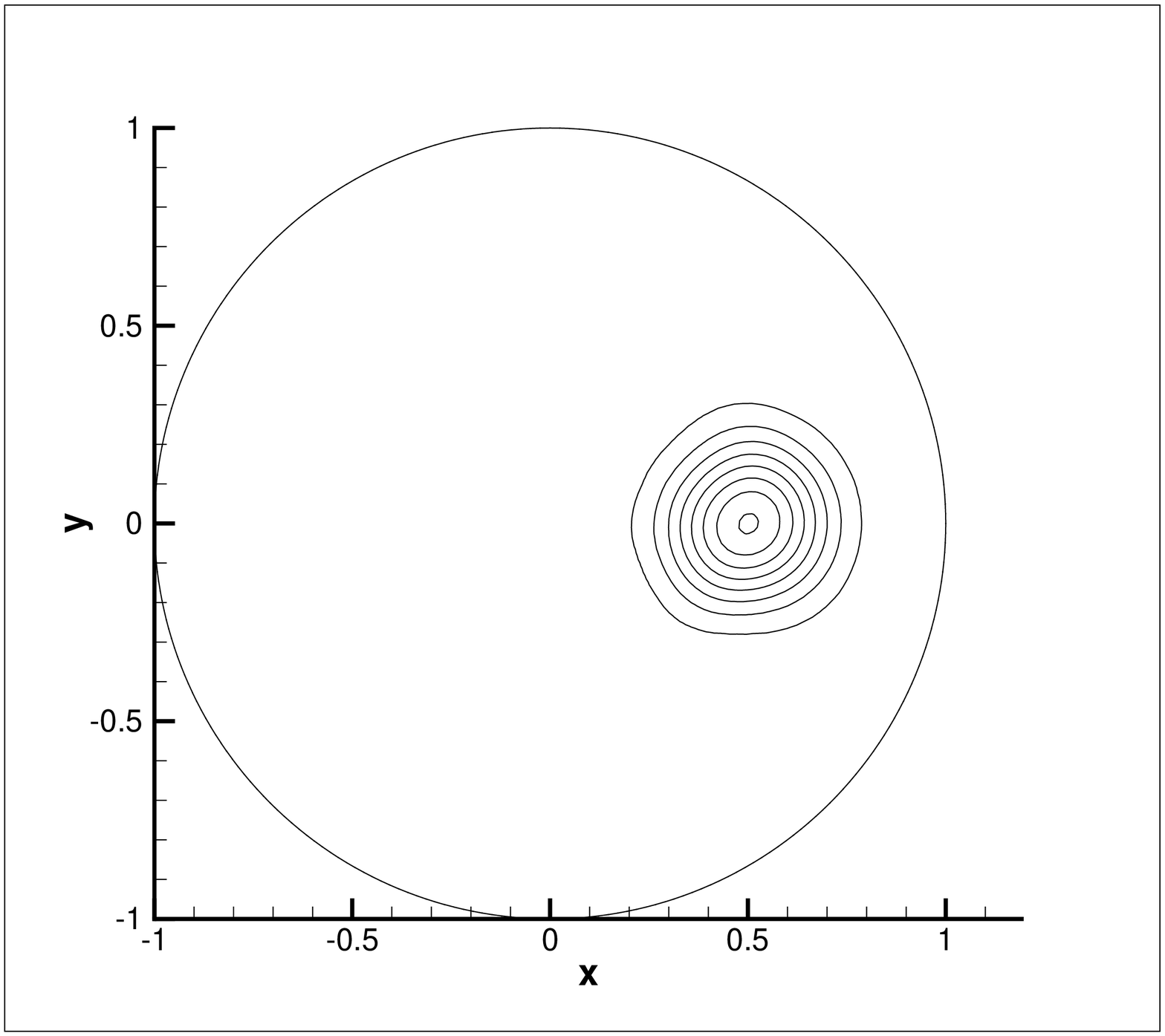}}\\
\multicolumn{2}{c}{(c)}
\end{tabular}
\end{center}
\caption{\label{abgrall_pleinary_P1} Results for the $\P^1$ approximation: (a) with SUPG, after 1 rotation, (b) with SUPG after 2 rotations,
 (c) with Galerkin+Jump after 10 rotations. The same isolines are represented.}
\end{figure}

\begin{figure}[h]
\begin{center}
\begin{tabular}{cc}
\includegraphics[width=0.45\textwidth]{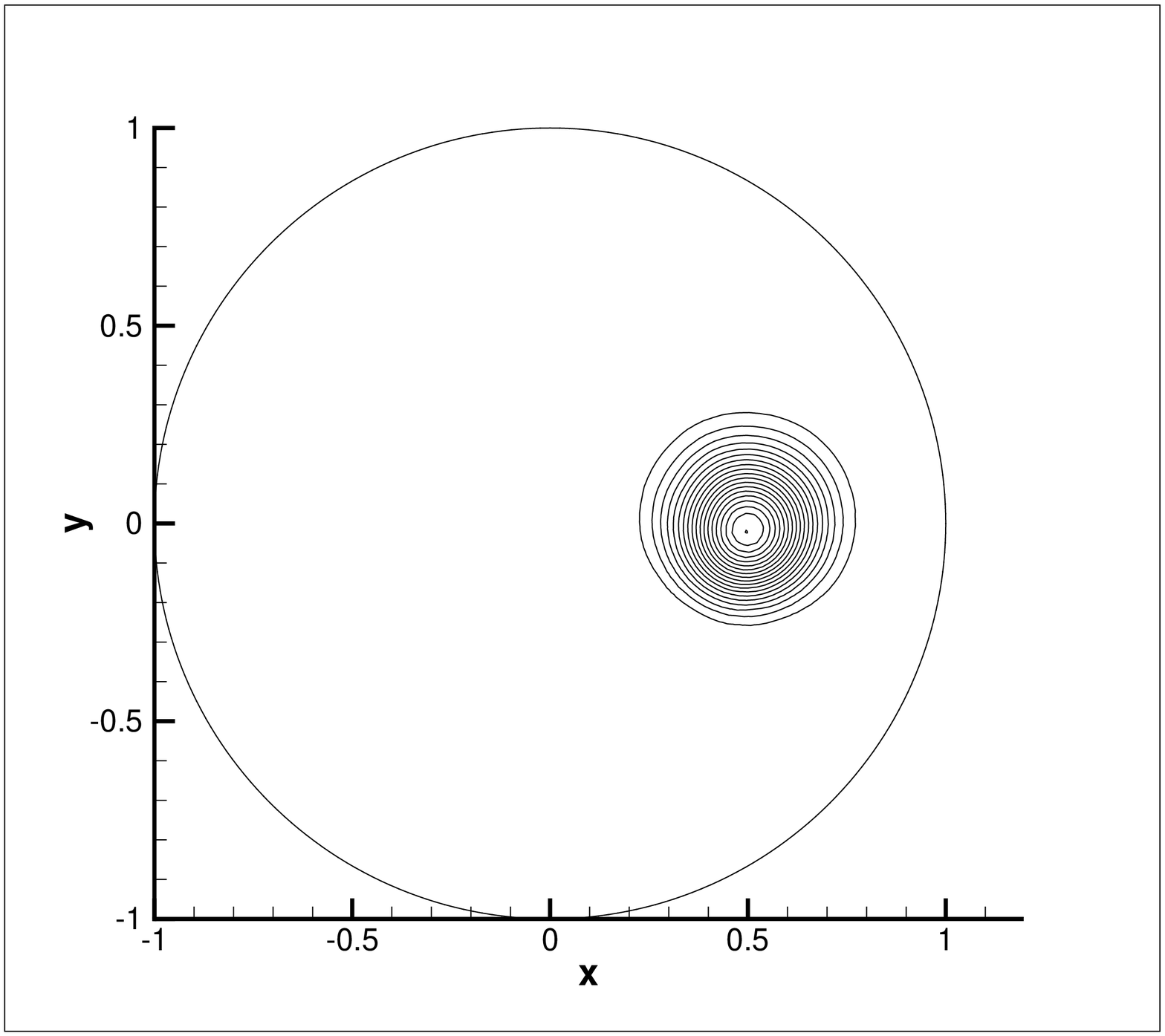}&
\includegraphics[width=0.45\textwidth]{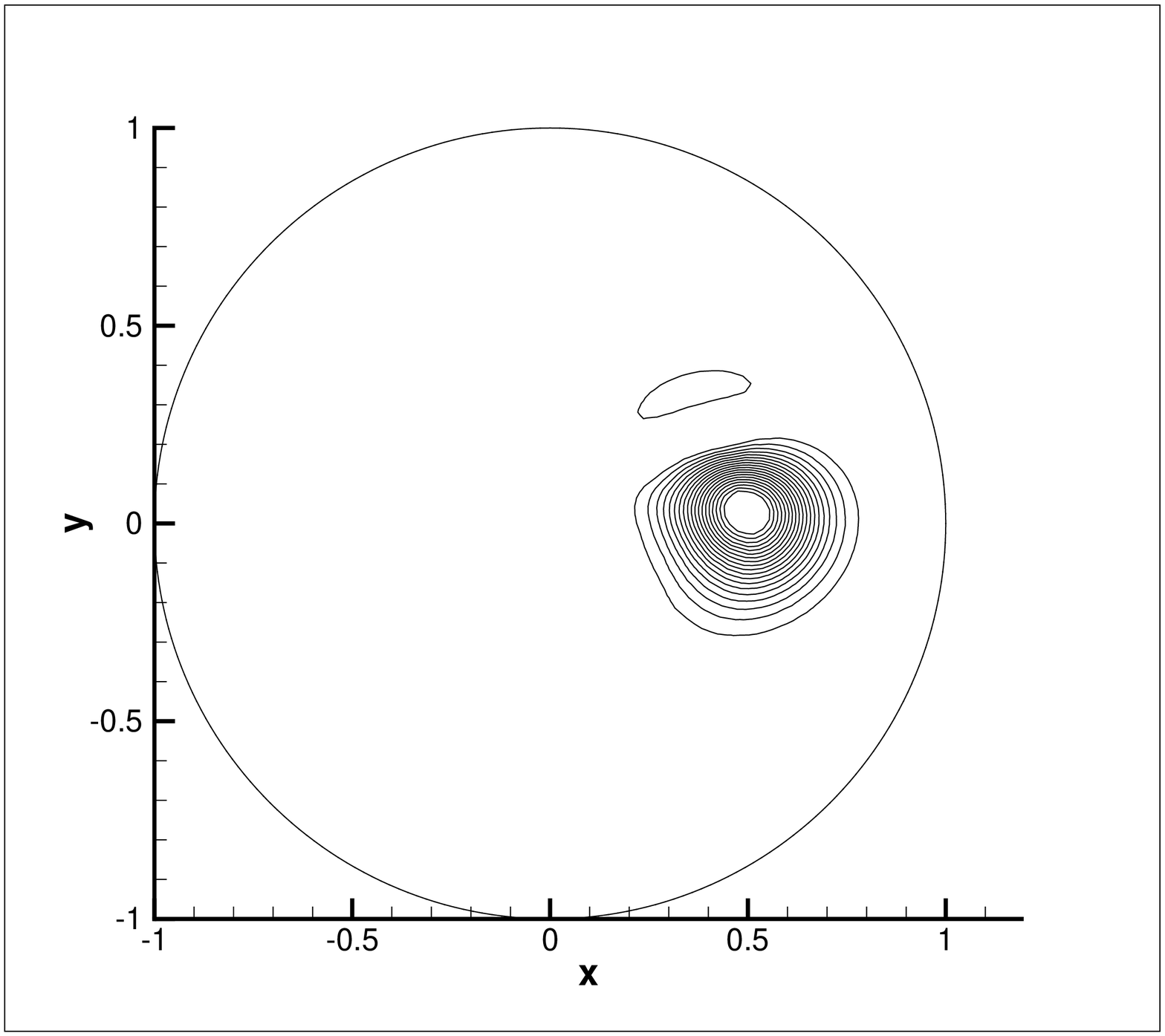}\\
(a)&(b)
\end{tabular}
\end{center}
\caption{\label{abgrall_pleinary_B2} Results for the $\mathbb{B}^2$ approximation: (a) with Galerkin+Jump after 10 rotations, (b) with SUPG after 1 rotations, (c). 
The same isolines are represented}
\end{figure}

Figure \ref{abgrall_pleinary_P1} displays the results for the $\P^1$ approximation, while Figure \ref{abgrall_pleinary_B2} shows those obtained for the quadratic approximation.
 The baseline schemes are the SUPG and the Galerkin scheme with jumps.

In Figure \ref{abgrall_pleinary_P1}, the same isolines are represented for the three results. We can see that after $10$ rotations, the results of the Galerkin+jump scheme look pretty good despite the coarse resolution.
 The minimum and maximum are $-0.012$ and $0.762$. For the SUPG results, after $1$ rotation, the minimum/maximum are $-0.004$ and $1.02$. After $2$ rotations we have $-0.047$ and
 $1.02$. This is better that what is obtained for Figure \ref{abgrall_pleinary_P1}-(c), but the dispersive effects are much more important for the SUPG scheme as it can be seen on Figure
 \ref{abgrall_pleinary_P1}-(b): this is why we have not shown further results for the SUPG/P1 case.

In Figure \ref{abgrall_pleinary_B2}, we show similar results obtained with the quadratic approximation. Again, the Galerkin+jump method is way less dispersive that the SUPG (stopped after only one rotation this time). We have found that if we perform $4$, $6$ or $8$ iterations of the defect correction, the quality of the SUPG improves a lot,
 but the cost becomes prohibitive with respect to the Galerkin+jump method for which, after $10$ rotations, the min/max are $-0.0044$ and $0.95$. We also see that the 
solution improves a lot with respect to linear elements, for example in terms of min/max values. There is however some dispersion, if we compare with the exact solution.

\section{Conclusions, perspectives}
The paper deals with the numerical approximation of linear scalar hyperbolic problems. We have shown, by carefully choosing the spatial approximation, 
and by using a non standard time step discretization, that it is possible to avoid the use of mass matrix in this problem, contrarily to what is usually thought about.
 The cost, on paper, is similar to a standard Runge-Kutta scheme, at least if we consider second and third order in time.  In a preliminary work, we have had
 similar results for the 1D advection problem, which are not shown here. We had also obtained the expected convergence slope.

A lot remains to be done. First, we have found experimental CFL conditions but this has to be rationalized by a numerical analysis. This method needs to
 be extended to non-linear problems. Preliminary results seems promising, but the results need to be checked on a wider range of problems, this is why we have
 not reported them here.
Last, this method needs to be extended to systems, for example the Euler equations of fluid mechanics.

\section*{Acknowledgements.} This research was funded by SNFS grant \# 200021\_15360. Early discussions with Mario Ricchiuto (INRIA Bordeaux Sud Ouest, France) are aknowledged.


\begin{thebibliography}{10}

\bibitem{abgrall_pleinary_abgrall:high}
R.~Abgrall.
\newblock Some comments about high order approximation of unsteady linear and
  non linear hyperbolic problems by continuous finite elements.
\newblock {\em In preparation}, 2016.

\bibitem{abgrall_pleinary_jirka}
R\'emi Abgrall and Jirka Trefilick.
\newblock An example of high order residual distribution scheme using
  non-Lagrange elements.
\newblock {\em Journal of Scientific Computing}, 45(1-3):64--89, October 2010.

\bibitem{abgrall_pleinary_burman}
E.~Burman and P.~Hansbo.
\newblock Edge stabilization for {G}alerkin approximation of
  convection-diffusin-reaction problems.
\newblock {\em Comput. Methods Appl. Mech. Engrg}, 193:1437--1453, 2004.

\bibitem{abgrall_pleinary_isoHughes}
J.~Austin Cottrell, Thomas~J.R. Hughes, and Yuri Bazilevs.
\newblock {\em Isogeometric Analysis: Toward Integration of CAD and FEA}.
\newblock John Wiley \& Sons, 2009.
\newblock ISBN 978-0-470-74873-2.

\bibitem{abgrall_pleinary_greengard}
A.~Dutt, L.~Greengard, and V.~Rokhlin.
\newblock Spectral deferred correction methods for ordinary differential
  equations.
\newblock {\em BIT Numerical Mathematics}, 40(2):241--266, 2000.

\bibitem{abgrall_pleinary_ho_lumping1}
G.Cohen, P.Joly, J.E.Roberts, and N.Tordjman.
\newblock High order triangular finite elements with mass lumping for the wave
  equation.
\newblock {\em SIAM J. Numer. Anal.}, 38(6):2047--2078, 2001.

\bibitem{abgrall_pleinary_Hughes1987}
T.J.R. Hughes and M.~Mallet.
\newblock A new finite element formulation for computational fluid dynamics
  {III}. The generalized streamline operator for multi-dimensional
  advective-diffusive systems.
\newblock {\em Comput. Methods Appl. Mech. Engrg.}, 58:305--328, 1987.

\bibitem{abgrall_pleinary_JohnsonNavertPitaka}
C.~Johnson, U.~N\"avert, and J.~Pitk\"aranta.
\newblock Finite element methods for linear hyperbolic problems.
\newblock {\em Comput. Methods Appl. Mech. Engrg}, 45:285--312, 1984.

\bibitem{abgrall_pleinary_ho_lumping2}
S.~Jund and S.~Salmon.
\newblock Arbitrary high order finite element schemes and high order mass
  lumping.
\newblock {\em Int. J. Appl. Math. Comput. Sci.}, 17(3):375--393, 2007.

\bibitem{abgrall_pleinary_mario}
Mario Ricchiuto and R\'emi Abgrall.
\newblock Explicit runge-kutta residual-distribution schemes for time dependent
  problems.
\newblock {\em Journal of Computational Physics}, 229(16):5653--5691, 2010.

\end{thebibliography}
\end{document}